\documentclass[leqno,draft]{article}

\def\supp{\mathop{\rm supp}}
\def\re{\mathop{\rm Re}}

\newtheorem{theorem}{Theorem}
\newtheorem{lemma}[theorem]{Lemma}
\newtheorem{proposition}[theorem]{Proposition}
\newtheorem{definition}[theorem]{Definition}
\newtheorem{corollary}[theorem]{Corollary}

\newcommand{\begintheorem}{\addtocounter{equation}{1}\begin{theorem}}
\newcommand{\beginlemma}{\addtocounter{equation}{1}\begin{lemma}}
\newcommand{\beginproposition}{\addtocounter{equation}{1}\begin{proposition}}
\newcommand{\begindefinition}{\addtocounter{equation}{1}\begin{definition}}
\newcommand{\begincorollary}{\addtocounter{equation}{1}\begin{corollary}}

\begin{document}

\title{Elements of harmonic analysis}

\author{Stephen William Semmes  \\
	Rice University		\\
	Houston, Texas}

\date{}

\maketitle

\tableofcontents

\bigskip

	These informal notes are based on a course given at Rice
University in the spring semester of 2004, and much more information
can be found in the references.

\section{Finite abelian groups}
\label{finite abelian groups}
\setcounter{equation}{0}

	Let $A$ be a finite abelian group.  Thus $A$ is a nonempty
set equipped with a binary operation which we denote $+$, which
is to say that if $a$, $b$ are elements of $A$, then $a + b$ is
a well-defined element of $A$.  This operation is commutative and
associative, which is to say that
\begin{equation}
	a + b = b + a
\end{equation}
and
\begin{equation}
	(a + b) + c = a + (b + c)
\end{equation}
for all $a, b, c \in A$.  There is an identity element $0$ in $A$,
which is characterized by the property
\begin{equation}
	a + 0 = a
\end{equation}
for all $a \in A$, and for each $a \in A$ there is a unique inverse
element, denoted $-a$, which is characterized by
\begin{equation}
	a + (-a) = 0.
\end{equation}

	By a \emph{character} on $A$ we mean a function $\phi$
from $A$ to the complex numbers with modulus $1$ such that
\begin{equation}
	\phi(a + b) = \phi(a) \, \phi(b)
\end{equation}
for all $a, b \in A$.  In other words, $\phi$ is a homomorphism from
$A$ into the group ${\bf T}$ of complex numbers with modulus $1$ using
multiplication as the group operation.  Let us note that if $\phi$ is
a homomorphism from $A$ into the nonzero complex numbers using
multiplication as the group operation, then $\phi$ automatically takes
values in ${\bf T}$.  Indeed, for each element $a$ of $A$, there is a
positive integer $m$ so that the sum of $a$ $n$ times is equal to $0$,
and this leads to the conclusion that $\phi(a)^m = 1$.  In other
words, $\phi(a)$ is an $m$th root of unity, and it has modulus equal
to $1$ in particular.

	Let $A^*$ denote the set of characters on $A$.  Note that we
automatically have a \emph{unit character} on $A$, which sends every
element of $A$ to $1$.  Also, if $\phi_1$, $\phi_2$ are characters on
$A$, then so is the product $\phi_1 \, \phi_2$, and if $\phi$ is a
character on $A$, then $1/\phi = \overline{\phi}$ is a character on
$A$.  Here we use the standard notation that if $z = x + i \, y$ is a
complex number with $x$, $y$ real numbers, the real and imaginary
parts of $z$, then $\overline{z}$ is the \emph{complex conjugate} of
$z$, given by $\overline{z} = x - i \, y$.  In short $A^*$ becomes
an abelian group with respect to multiplication of characters, called
the \emph{dual} of $A$.

	Let $V$ denote the vector space of complex-valued functions
on $A$.  For each $a \in A$ we can define a linear transformation
$T_a$ on $V$ by
\begin{equation}
	T_a(f)(x) = f(x - a)
\end{equation}
for each $f \in V$.  In other words, $T_a$ translates a given function
by $a$.  Notice that
\begin{equation}
	T_a \circ T_b = T_{a + b}
\end{equation}
for all $a, b \in A$, and that for each $a \in A$ $T_a$ is an
invertible linear operator on $V$, with 
\begin{equation}
	(T_a)^{-1} = T_{-a}.
\end{equation}

	If $\phi \in A^*$, then
\begin{equation}
	T_a(\phi) = \overline{\phi(a)} \, \phi,
\end{equation}
since $\phi(x-a) = \phi(-a) \, \phi(x)$ and $\phi(-a) = \phi(a)^{-1} =
\overline{\phi(a)}$.  Thus $\phi$ is an eigenvector for $T_a$ with
eigenvalue $\overline{\phi(a)}$.  Now suppose that $f(x)$ is a
function on $A$ which is an eigenvector for $T_a$ for each $a \in A$
with eigenvalue $\lambda(a)$, which is to say that
\begin{equation}
	T_a(f) = \lambda(a) \, f
\end{equation}
for all $a \in A$.  This is equivalent to
\begin{equation}
	f(x - a) = \lambda(a) \, f(x)
\end{equation}
for all $x, a \in A$.  It is easy to see that $f$ is either identically
$0$ on $A$, or that $\lambda$ defines a character on $A$ and $f$
is a nonzero multiple of the complex conjugate of $\lambda$.

	Let $n$ denote the number of elements of $A$.  If $f_1$, $f_2$
are elements of $A$, then let us define their inner product
$\langle f_1, f_2 \rangle_A$ by
\begin{equation}
	\langle f_1, f_2 \rangle_A 
		= \frac{1}{n} \sum_{x \in A} f_1(x) \, \overline{f_2(x)}.
\end{equation}
If $f$ is a function on $V$, then we also define its norm $\|f\|_A$ by
\begin{equation}
	\|f\|_A = \biggl(\frac{1}{n} \sum_{x \in A} |f(x)|^2 \biggr)^{1/2},
\end{equation}
which is the same as
\begin{equation}
	\|f\|_A = \langle f, f \rangle_A^{1/2}.
\end{equation}
By standard results we have the \emph{Cauchy--Schwarz inequality}
\begin{equation}
	|\langle f_1, f_2 \rangle_A| \le \|f_1\|_A \, \|f_2\|_A
\end{equation}
and the \emph{triangle inequality}
\begin{equation}
	\|f_1 + f_2\|_A \le \|f_1\|_A + \|f_2\|_A
\end{equation}
for $f_1, f_2 \in A$.

	With respect to this inner product, each translation operator
$T_a$ is a \emph{unitary} operator on $V$, so that
\begin{equation}
	\langle T_a(f_1), T_a(f_2) \rangle_A
		= \langle f_1, f_2 \rangle_A
\end{equation}
for all $f_1, f_2 \in A$.  It is a well-known result from linear
algebra that a unitary linear operator on a finite-dimensional
complex inner product space can be diagonalized in an orthonormal
basis.  Moreover, given a finite collection of unitary transformations
which commute with each other, there is in fact an orthonormal
basis for the inner product space in which all of the unitary
transformations are diagonalized.  Of course the eigenvalues
of a unitary transformation automatically have modulus equal to $1$.

	Suppose that $\phi$ is a character on $A$.  If $a$ is
any element of $A$, then
\begin{equation}
	\phi(a) \sum_{x \in A} \phi(x) = \sum_{x \in A} \phi(x + a)
		= \sum_{x \in A} \phi(x),
\end{equation}
where the second equality holds because the two sums have the same
terms, just arranged differently.  It follows that either $\phi(a) =
1$ for all $a \in A$, so that $\phi$ is the unit character, or that
\begin{equation}
	\sum_{x \in A} \phi(x) = 0.
\end{equation}
If $\psi_1, \psi_2 \in A^*$, then we can apply this to the character
\begin{equation}
	\phi(x) = \psi_1(x) \, \overline{\psi_2(x)}
\end{equation}
to conclude that either $\psi_1 = \psi_2$ or
\begin{equation}
	\langle \psi_1, \psi_2 \rangle_A = 0,
\end{equation}
which is to say that $\psi_1$, $\psi_2$ are orthogonal to each other.
Of course this also follows from the characterization of characters as
simultaneous eigenvectors of the translation operators, since distinct
characters correspond to distinct eigenvalues for at least one
eigenvector.  Because of the way that we defined the inner product
on $V$, characters automatically have norm equal to $1$.

\begintheorem
The dual group $A^*$ has the same number of elements as $A$ does,
and the elements of $A^*$ form an orthonormal basis for $V$.
\end{theorem}

	More precisely, we have seen that the characters are
orthonormal with respect to the inner product $\langle \cdot, \cdot
\rangle_A$, and hence that they are linearly independent, which
implies that there are at most $n$ elements of $A^*$.  We also
know that simultaneous eigenvectors of the translation operators
are multiples of characters, and since there is an orthonormal
basis of simultaneous eigenvectors, we get an orthonormal basis
of characters by multiplying the simultaneous eigenvectors
in the orthonormal basis by complex numbers with modulus $1$,
if necessary.  In particular, the number of characters is exactly
equal to $n$.

	Let $W$ denote the vector space of complex-valued functions on
$A^*$, which therefore has the same dimension as $V$, $n$.  If $h_1$,
$h_2$ are elements of $W$, let us define their inner product to be
\begin{equation}
	\langle h_1, h_2 \rangle_{A^*}
		= \sum_{\phi \in A^*} h_1(\phi) \overline{h_2(\phi)}.
\end{equation}
This leads to the norm
\begin{equation}
	\|h\|_{A^*} = \langle h, h \rangle_{A^*}^{1/2}
		= \biggl(\sum_{\phi \in A^*} |h(\phi)|^2 \biggr)^{1/2}.
\end{equation}
As before we have the Cauchy-Schwarz and triangle inequalities, which
is to say that
\begin{equation}
	|\langle h_1, h_2 \rangle_{A^*}|
		\le \|h_1\|_{A^*} \, \|h_2\|_{A^*}
\end{equation}
and
\begin{equation}
	\|h_1 + h_2\|_{A^*} \le \|h_1\|_{A^*} + \|h_2\|_{A^*}
\end{equation}
for all $h_1, h_2 \in W$.

	If $f$ is a function on $A$, then the \emph{Fourier transform}
of $f$ is the function $\widehat{f}(\phi)$ on $A^*$ defined by
\begin{equation}
	\widehat{f}(\phi) = \langle f, \phi \rangle_A
\end{equation}
for $\phi \in A$.  We may also denote the Fourier transform of $f$ by
$\mathcal{F}(f)$, so that $\mathcal{F}$ is a linear transformation
from $V$ to $W$.  The orthonormality of the characters on $A$ implies
that the Fourier transform is actually a unitary transformation from
$V$ to $W$ with respect to the inner products that we have defined,
i.e.,
\begin{equation}
	\langle \mathcal{F}(f_1), \mathcal{F}(f_2) \rangle_{A^*}
		= \langle f_1, f_2 \rangle_A
\end{equation}
for all $f_1, f_2 \in V$, and thus
\begin{equation}
	\|\mathcal{F}(f)\|_{A^*} = \|f\|_A
\end{equation}
for all $f \in V$.

	As a Fourier inversion formula we can write
\begin{equation}
	f = \sum_{\phi \in A^*} \widehat{f}(\phi) \, \phi
\end{equation}
for all $f \in V$.  In particular a function $f$ on $A$ is uniquely
determined by its Fourier transform, which is to say that the Fourier
transform is a one-to-one linear mapping from $V$ into $W$.  It
follows that the Fourier transform maps $V$ onto $W$, since they have
the same dimension.  In other words, every function on $A^*$ arises as
the Fourier transform of a function on $A$, which can be obtained
simply by using the function on $A^*$ as the coefficients for an
expansion in characters of a function on $A$.

	If $f_1, f_2 \in V$, then the \emph{convolution} of $f_1$
and $f_2$ is the function on $A$ defined by
\begin{equation}
	(f_1 * f_2)(x) = \frac{1}{n} \sum_{x \in A} f_1(y) \, f_2(x - y).
\end{equation}
This is clearly linear in each of $f_1$ and $f_2$.  It is also
commutative and associative, which is to say that
\begin{equation}
	f_1 * f_2 = f_2 * f_1
\end{equation}
and
\begin{equation}
	(f_1 * f_2) * f_3 = f_1 * (f_2 * f_3)
\end{equation}
for all $f_1, f_2, f_3 \in V$.  Moreover,
\begin{equation}
	\mathcal{F}(f_1 * f_2) = \mathcal{F}(f_1) \, \mathcal{F}(f_2)
\end{equation}
for all $f_1, f_2 \in V$, i.e., the Fourier transform of a convolution
is the same as the product of the individual Fourier transforms.

\section{Riesz--Thorin convexity}
\label{Riesz--Thorin convexity}
\setcounter{equation}{0}

	Fix a positive integer $m$, and consider ${\bf C}^m$
as an $m$-dimensional vector space over the complex numbers.
If $p$ is a real number such that $1 \le p < \infty$ and
$v = (v_1, \ldots, v_m)$ is a vector in ${\bf C}^m$, put
\begin{equation}
	\|v\|_p = \biggl(\sum_{j=1}^m |v_j|^p \biggr)^{1/p}.
\end{equation}
When $p = \infty$ set
\begin{equation}
	\|v\|_\infty = \max \{|v_j| : 1 \le j \le m\}.
\end{equation}

	For $1 \le p \le \infty$ and $v \in {\bf C}^m$ we have that
$\|v\|_p$ is a nonnegative real number which is equal to $0$ if and
only if $v = 0$, and that
\begin{equation}
	\|\alpha \, v\|_p = |\alpha| \, \|v\|_p
\end{equation}
for all complex numbers $\alpha$.  The triangle inequality
\begin{equation}
	\|v + w\|_p \le \|v\|_p + \|w\|_p
\end{equation}
for all $v, w \in {\bf C}^m$ is equivalent to the convexity of the
closed unit ball
\begin{equation}
	B_p = \{v \in {\bf C}^m : \|v\|_p \le 1\},
\end{equation}
which is to say that $t \, v + (1 - t) \, w \in B_p$ whenever
$v, w \in B_p$ and $0 \le t \le 1$.  When $p = 1, \infty$ it
is quite easy to verify the triangle inequality directly,
while for $1 < p < \infty$ one can derive it from the
convexity of the function $x^p$ on $[0, \infty)$.

	It is clear from the definitions that
\begin{equation}
	\|v\|_\infty \le \|v\|_p
\end{equation}
when $v \in {\bf C}^m$ and $1 \le p < \infty$.  Using this one
can check that
\begin{equation}
	\|v\|_q \le \|v\|_p
\end{equation}
when $v \in {\bf C}^m$ and $1 \le p \le q < \infty$.  Namely,
\begin{eqnarray}
	\|v\|_q^q = \sum_{j=1}^m |v_j|^q
		& \le & \|v\|_\infty^{q-p} \sum_{j=1}^m |v_j|^p	\\
		& \le & \|v\|_p^{q-p} \, \|v\|_p^p = \|v\|_p^p.
						\nonumber
\end{eqnarray}

	In the other direction,
\begin{equation}
	\|v\|_p \le m^{1/p} \, \|v\|_\infty
\end{equation}
when $v \in {\bf C}^m$ and $1 \le p < \infty$.  This follows easily
from the definitions.  A slightly more tricky fact is that
\begin{equation}
	\|v\|_p \le m^{1/p - 1/q} \, \|v\|_\infty
\end{equation}
for $1 \le p \le q < \infty$.  This can be derived from the convexity
of the function $x^{q/p}$ on $[0, \infty)$.

	If $1 \le p \le \infty$, then the \emph{conjugate exponent}
$p'$, $1 \le p' \le \infty$, is defined by the condition
\begin{equation}
	\frac{1}{p} + \frac{1}{p'} = 1.
\end{equation}
\emph{H\"older's inequality} states that
\begin{equation}
	\biggl| \sum_{j=1}^m v_j \, w_j \biggr|
		\le \|v\|_p \, \|w\|_{p'}
\end{equation}
for all $v, w \in {\bf C}^m$.  When $p = 1$, $p' = \infty$
or $p = \infty$, $p' = 1$ this is a straightforward consequence
of the definitions, and when $p = p' = 2$ this is the classical
Cauchy--Schwarz inequality.

	Suppose that $1 < p < \infty$, so that $1 < p' < \infty$, and
let us prove H\"older's inequality.  If $x$, $y$ are nonnegative real
numbers, then we have that
\begin{equation}
	x \, y \le \frac{x^p}{p} + \frac{y^{p'}}{p'}.
\end{equation}
This can be derived from the convexity of the exponential function,
for instance.  For $v, w \in {\bf C}^m$ it follows that
\begin{equation}
	\biggl|\sum_{j=1}^m v_j \, w_j \biggr|
		\le \sum_{j=1}^m |v_j| \, |w_j|
		\le \frac{\|v\|_p^p}{p} + \frac{\|w\|_{p'}}{p'},
\end{equation}
by applying the previous inequality to each of the terms in the sum.
This implies H\"older's inequality when $\|v\|_p$, $\|w\|_{p'}$
are equal to $1$, and one can reduce to that case using the
homogeneity of the norms.

	H\"older's inequality is sharp, in the sense that for each $v
\in {\bf C}^m$ and $1 \le p \le \infty$ there is a $w \in {\bf C}^m$
such that $\|w\|_{p'} = 1$ and
\begin{equation}
	\sum_{j=1}^m v_j \, w_j = \|v\|_p.
\end{equation}
This is easy to verify from the definitions.  Notice that if $1 < p <
\infty$ and $v \ne 0$, then $w$ is unique.

	Now let $T$ be a linear transformation from ${\bf C}^m$
into itself.  Suppose that
\begin{equation}
	1 \le p_0, p_1, q_0, q_1 \le \infty,
\end{equation}
and let $p_0'$, $p_1'$, $q_0'$, $q_1'$ be their conjugate exponents.
Suppose also that $L_0$, $L_1$ are nonnegative real numbers such that
\begin{equation}
	\|T(v)\|_{q_0} \le L_0 \, \|v\|_{p_0}
\end{equation}
and
\begin{equation}
	\|T(v)\|_{q_1} \le L_1 \, \|v\|_{p_1}
\end{equation}
for all $v \in {\bf C}^m$.  The \emph{Riesz--Thorin convexity
theorem} states that
\begin{equation}
	\|T(v)\|_{q_t} \le L_t \, \|v\|_{p_t}
\end{equation}
for $0 < t < 1$, where $1/p_t = (1-t)/ p_0 + t/p_1$, $1/q_t =
(1-t)/q_0 + t/q_1$, and $L_t = L_0^{(1-t)} \, L_1^t$.

	Here is an equivalent formulation, using H\"older's inequality
and the fact that equality is attained in H\"older's inequality.
Suppose that
\begin{equation}
	\biggl|\sum_{j=1}^m (T(v))_j \, w_j \biggr|
		\le L_0 \, \|v\|_{p_0} \, \|w\|_{q_0'}
\end{equation}
and
\begin{equation}
	\biggl|\sum_{j=1}^m (T(v))_j \, w_j \biggr|
		\le L_1 \, \|v\|_{p_1} \, \|w\|_{q_1'}
\end{equation}
for all $v, w \in {\bf C}^m$, where $(T(v))_j$ denotes the
$j$th component of $T(v)$.  Then
\begin{equation}
	\biggl|\sum_{j=1}^m (T(v))_j \, w_j \biggr|
		\le L_t \, \|v\|_{p_t} \, \|w\|_{q_t'}
\end{equation}
for all $v, w \in {\bf C}^m$ and $0 < t < 1$, with $q_t'$ being
the conjugate exponent of $q_t$.

	Marcel Riesz proved this originally using a real-variable
method that worked for both ${\bf R}^m$ and ${\bf C}^m$, and with the
additional hypothesis that $q_0 \ge p_0$, $q_1 \ge p_1$.  Thorin then
found a way to use complex analytic functions to prove this in a nice
way and without this extra condition on the exponents, in the complex
case.  In general the real case can be reduced to the complex case
with a slightly less sharp inequality.  We shall now discuss Thorin's
argument.

	We begin with some preliminary facts from complex analysis.
Let $U$ be a nonempty bounded open subset of the complex plane, and
let $f(z)$ be a continuous complex-valued function on the closure
$\overline{U}$ of $U$ which is complex-analytic on $U$.  Because
$U$ is bounded, $\overline{U}$ and $\partial U$ are closed and bounded
subsets of the complex plane which are therefore compact, and
the maximum of $|f(z)|$ on $\overline{U}$ and $\partial U$ are
attained.  The \emph{maximum principle} for complex analytic
functions implies that
\begin{equation}
	\max \{|f(z)| : z \in \overline{U} \}
		\le \max \{|f(z)| : z \in \partial U\}.
\end{equation}

	We would like to apply this to the case where $U$ is the strip
$\{z \in {\bf C} : 0 < \re z < 1\}$, where $\re z$ denotes the real
part of a complex number $z$.  This region is not bounded, and we
shall deal with that in a moment.  Let $f(z)$ be a continuous
complex-valued function defined on the closure of $U$, which is the
set of $z \in {\bf C}$ such that $0 \le \re z \le 1$, and suppose that
$f(z)$ is holomorphic on $U$ and bounded on the closure of $U$.
For $0 \le t \le 1$ put
\begin{equation}
	M_t = \sup \{|f(z)| : \re z = t\}.
\end{equation}
We would like to say that
\begin{equation}
	M_t \le \max (M_0, M_1)
\end{equation}
when $0 < t < 1$, which would say that the supremum of $|f(z)|$
over the closure of $U$ is equal to the supremum of $|f(z)|$
over the boundary of $U$.

	If $|f(z)| \to 0$ as $|z| \to \infty$, then this extension of
the maximum principle can be derived from the version for bounded
regions, by approximating the strip $U$ with large rectangular
subsets.  For a general bounded function $f(z)$ we can reduce
to this case using the following trick.  For each $\epsilon > 0$,
consider the function
\begin{equation}
	f_\epsilon(z) = \exp (\epsilon \, z^2) \, f(z).
\end{equation}
This is also a continuous complex-valued function on the closure
of $U$ which is complex-analytic on $U$.  If $z = x + i \, y$,
with $x, y \in {\bf R}$, then
\begin{equation}
	|\exp (\epsilon \, z^2)| = \exp (\epsilon (x^2 - y^2)),
\end{equation}
and one can derive the maximum principle for $f$ on the closure of $U$
from the corresponding statement for $f_\epsilon$ for all $\epsilon >
0$.

	As a refinement of the maximum principle in this situation
there is the \emph{three lines theorem}, which states that if $f(z)$
is a bounded continuous complex-analytic function on the closure of
$U$ which is complex-analytic on $U$, and if $M_t$ is as defined
previously, then
\begin{equation}
	M_t \le M_0^{(1-t)} \, M_1^t
\end{equation}
when $0 < t < 1$.  This can be derived by applying the previous
version to functions of the form $\exp (a \, z) \, f(z)$,
where $a$ is a real number.  The main point is that
\begin{equation}
	|\exp (a \, z)| = \exp (a \, x)
\end{equation}
if $x = \re z$, and in particular this quantity is constant
on the lines $\re z = t$.

	Now let us return to the Riesz--Thorin convexity theorem.  Let
$v, w \in {\bf C}^m$ and $0 < t < 1$ be given.  Without loss of
generality let us make the normalizing assumption that
\begin{equation}
	\|v\|_{p_t} = \|w\|_{q_t'} = 1.
\end{equation}
We would like to show that
\begin{equation}
	\biggl| \sum_{j=1}^m (T(v))_j \, w_j \biggr|
		\le L_0^{(1-t)} \, L_1^t.
\end{equation}

	The idea is to realize this as a case of the three lines
theorem.  To do this we would like to find bounded continuous ${\bf
C}^m$-valued functions $\alpha(z)$, $\beta(z)$ defined on the closed
unit strip $\{z \in {\bf C} : 0 \le \re z \le 1 \}$ which are
complex-analytic in the interior and satisfy
\begin{equation}
	\alpha(t) = v, \quad \beta(t) = w
\end{equation}
and
\begin{equation}
	\|\alpha(z)\|_{p_x} = \|\beta(z)\|_{q_x'} = 1
\end{equation}
when $0 \le x \le 1$ and $\re z = x$.  If we can do this, then the
desired inequality follows from the three lines theorem applied to
\begin{equation}
	\sum_{j=1}^m (T(\alpha(z)))_j \, \beta_j(z).
\end{equation}
Here and in the following we write $\alpha_j(z)$, $\beta_j(z)$
for the components of $\alpha(z)$, $\beta(z)$.

	When $v_j = 0$ or $w_k = 0$ we put $\alpha_j(z) = 0$
and $\beta_k(z) = 0$ for all $z$.  Otherwise we put
\begin{equation}
	\alpha_j(z) = v_j \, |v_j|^{a_0 \, z + a_1}, \quad
		\beta_k(z) = w_k \, |w_k|^{b_0 \, z + b_1},
\end{equation}
where $a_0$, $a_1$, $b_0$, $b_1$ are fixed real numbers.  More
precisely, they should be chosen so that
\begin{equation}
	\frac{1}{p_t} (a_0 \, x + a_1 + 1) = \frac{1}{p_x},
   \quad \frac{1}{q_t'} (b_0 \, x + b_1 + 1) = \frac{1}{q_x'}
\end{equation}
when $0 \le x \le 1$.  The details are left as an exercise to the
interested reader.

	Now let us consider applications of the Riesz--Thorin theorem
to Fourier analysis on finite abelian groups as in the previous
section.  Let $A$ be a finite abelian group with $n$ elements, and let
$A^*$ be the dual group.  Also let $V$ denote the vector space of
complex-valued functions on $A$, and let $W$ denote the vector space
of complex-valued functions on $A^*$, each of which can be identified
with ${\bf C}^n$.

	If $f \in V$, put
\begin{equation}
	\|f\|_{p, A} = \biggl(\frac{1}{n} \sum_{x \in A} |f(x)|^p \biggr)^{1/p}
\end{equation}
when $1 \le p < \infty$ and
\begin{equation}
	\|f\|_{\infty, A} = \max \{|f(x)| : x \in A\}.
\end{equation}
Thus $\|f\|_{p, A}$ differs from the norm previously defined by a
constant factor, which does not really cause any trouble.  For
$h \in W$ we put
\begin{equation}
	\|h\|_{p, A^*} = \biggl(\sum_{\phi \in A^*} |h(\phi)|^p \biggr)^{1/p}
\end{equation}
when $1 \le p < \infty$ and
\begin{equation}
	\|h\|_{\infty, A^*} = \max \{|h(\phi)| : \phi \in A^*\},
\end{equation}
which is exactly the same as the norms defined earlier in this section.

	For the record, if $f \in V$ then
\begin{equation}
	\|f\|_{p, A} \le \|f\|_{q, A} \le n^{1/p - 1/q} \, \|f\|_{p, A}
\end{equation}
when $1 \le p \le q \le \infty$.  For $h \in W$ we still have
\begin{equation}
	\|h\|_{q, A^*} \le \|h\|_{p, A^*} \le n^{1/p - 1/q} \, \|h\|_{q, A^*}
\end{equation}
when $1 \le p \le q \le \infty$.  We can write H\"older's inequality
as
\begin{equation}
	\frac{1}{n} \biggl|\sum_{x \in A} f_1(x) \, f_2(x) \biggr|
		\le \|f_1\|_{p, A} \, \|f_2\|_{p', A}
\end{equation}
for all $f_1, f_2 \in A$ and
\begin{equation}
	\biggl|\sum_{\phi \in A^*} h_1(\phi) \, h_2(\phi) \biggr|
		\le \|h_1\|_{p, A^*} \, \|h_2\|_{p', A^*}
\end{equation}
for all $h_1, h_2 \in W$.  Here $1 \le p \le \infty$ and $p'$ is
the conjugate exponent of $p$.

	Let $\mathcal{F} : V \to W$ denote the Fourier transform.
Thus
\begin{equation}
	\|\mathcal{F}(f)\|_{2, A^*} = \|f\|_{2, A}
\end{equation}
for all $f \in V$, since when $p = 2$ the norms defined here reduce
to the ones discussed in the previous section, associated to the
inner products on $V$, $W$.  It is easy to see that
\begin{equation}
	\|\mathcal{F}(f)\|_{\infty, A^*} \le \|f\|_{1, A}
\end{equation}
for all $f \in V$, since the characters on $A$ have modulus equal
to $1$ at every point in $A$.  The convexity theorem implies that
\begin{equation}
	\|\mathcal{F}(f)\|_{p', A^*} \le \|f\|_{p, A}
\end{equation}
for all $f \in V$ when $1 \le p \le 2$.

	Now let us consider some inequalities related to convolutions.
Let $f_1$, $f_2$ be functions on $V$.  If $1 \le p \le \infty$, then
one can check that
\begin{equation}
	\|f_1 * f_2\|_{p, A} \le \|f_1\|_{1, A} \, \|f_2\|_{p, A}.
\end{equation}
Indeed, one can think of the convolution $f_1 * f_2$ as being a
linear combination of translates of $f_2$ with coefficients given
by the values of $f_1$.  The translated of $f_2$ have the same
$\|\cdot \|_{p, A}$ norm as $f_2$, so that the $\|\cdot \|_{p, A}$ norm
of $f_1 * f_2$ can be estimated as above using the homogeneity and
triangle inequality for the norm.

	Using H\"older's inequality it follows that
\begin{equation}
	\|f_1 * f_2\|_{\infty, A} \le \|f_1\|_{p', A} \, \|f_2\|_{p, A}.
\end{equation}
Let us think of $f_2$ as being fixed, and
\begin{equation}
	f_1 \mapsto f_1 * f_2
\end{equation}
as being a linear mapping from $V$ to itself.  Let us also fix $p$, so
that we have two estimates for this linear mapping as in the preceding
inequalities.  From the convexity theorem it follows that if $1 \le r
\le p'$ and $q$ is defined by
\begin{equation}
	\frac{1}{q} = \frac{1}{p} + \frac{1}{r} - 1,
\end{equation}
then
\begin{equation}
	\|f_1 * f_2 \|_{q, A} \le \|f_1\|_{r, A} \, \|f_2\|_{p, A}.
\end{equation}

	Fix a function $b$ on $A$, and define a linear transformation
$T_b : V \to V$ by
\begin{equation}
	T_b(f) = f * b.
\end{equation}
As above, if $1 \le p \le \infty$, then
\begin{equation}
	\|T_b(f)\|_{p, A} \le \|b\|_{1, A} \, \|f\|_{p, A}.
\end{equation}
When $p = 1$ this is sharp, because if $\delta_0$ denotes the function
on $A$ such that $\delta_0(x) = n$ when $x = 0$ and $\delta_0(x) = 0$
when $x \ne 0$, then
\begin{equation}
	T_b(\delta_0) = b
\end{equation}
and
\begin{equation}
	\|T_b\|_{1, A} = \|b\|_{1, A} = \|b\|_{1, A} \, \|\delta_0\|_{1, A}.
\end{equation}

	Now suppose that $p = 2$.  For this we use the Fourier
transform, the fact that it converts convolutions into
multiplications, and the unitary property for the norms.  Namely,
\begin{equation}
	\mathcal{F}(T_b(f)) = \mathcal{F}(b) \, \mathcal{F}(f)
\end{equation}
and
\begin{equation}
	\|\mathcal{F}(T_b(f))\|_{2, A^*}
		= \|T_b(f)\|_{2, A},
\end{equation}
and therefore
\begin{equation}
	\|T_b(f)\|_{2, A} \le \|\mathcal{F}(b)\|_{\infty, A^*} \, \|f\|_{2, A}.
\end{equation}
To put it another way, $T_b$ is diagonalized by the orthonormal basis
of $V$ consisting of the characters, with the diagonal entries given
by the values of the Fourier transform of $b$.  This inequality is
also sharp, with equality being attained for characters $\phi$ on $A$
such that
\begin{equation}
	|\mathcal{F}(\phi)|
\end{equation}
is as large as possible.

	We can apply the convexity theorem to these two estimates for
$p = 1$ and $p = 2$ to obtain that
\begin{equation}
	\|T_b(f)\|_{p, A}
		\le \|b\|_{1, A}^{2/p - 1} \, 
			\|\mathcal{F}(b)\|_{\infty, A^*}^{2 - 2/p}
				\, \|f\|_{p, A}
\end{equation}
when $1 \le p \le 2$.  When $p \ge 2$ the estimate is the same as for
$p'$, $1 \le p' \le 2$, by a duality argument.  To be more precise,
if $p$ is fixed and $M$ is a nonnegative real number such that
\begin{equation}
	\|T_b(f)\|_{p, A} \le M \, \|f\|_{p, A}
\end{equation}
for all $f \in V$, then the analogous inequality holds also for $p'$,
i.e.,
\begin{equation}
	\|T_b(f)\|_{p', A} \le M \, \|f\|_{p', A}
\end{equation}
for all $f \in V$.

	Indeed, if $f_1$, $f_2$ are functions on $A$, then
\begin{eqnarray}
	\frac{1}{n} \sum_{x \in A} T_b(f_1)(x) \, f_2(x)
  & = & 
 \frac{1}{n^2} \sum_{x \in A} \sum_{y \in A} b(y) \, f_1(x - y) \, f_2(x)
						\\
  & = & 
 \frac{1}{n^2} \sum_{x \in A} \sum_{y \in A} b(x - y) \, f_1(y) \, f_2(x)
					\nonumber \\
  & = & \frac{1}{n} \sum_{y \in A} f_1(y) \, T_{\widetilde{b}}(f_2)(y),
					\nonumber
\end{eqnarray}
where
\begin{equation}
	\widetilde{b}(z) = b(-z).
\end{equation}
Using this one can check that
\begin{equation}
	\|T_b(f_1)\|_{p, A} \le M \, \|f_1\|_{p, A}
\end{equation}
for all $f_1 \in V$ holds if and only if
\begin{equation}
	\|T_{\widetilde{b}}(f_2)\|_{p', A} \le M \, \|f_2\|_{p', A}
\end{equation}
for all $f_2 \in V$.  Also, for any $q$,
\begin{equation}
	\|T_{\widetilde{b}}(f)\|_{q, A} \le M \, \|f\|_{q, A}
\end{equation}
holds for all $f \in V$ if and only if
\begin{equation}
	\|T_b(f)\|_{q, A} \le M \, \|f\|_{q, A}.
\end{equation}
This follows from the fact that
\begin{equation}
	T_{\widetilde{b}}(f)(-x) = T_b(\widetilde{f})(x).
\end{equation}

\section{Continuous functions on ${\bf R}^n$}
\label{continuous functions on R^n}
\setcounter{equation}{0}

	Recall that if $x = (x_1, \ldots, x_n)$ is an element of ${\bf
R}^n$, then the standard Euclidean norm of $x$ is defined by
\begin{equation}
	|x| = \biggl(\sum_{j=1}^n |x_j|^2 \biggr)^{1/2}.
\end{equation}
If $x$, $y$ are elements of ${\bf R}^n$, then their inner product
is defined by
\begin{equation}
	\langle x, y \rangle = \sum_{j=1}^n x_j \, y_j,
\end{equation}
and we have that
\begin{equation}
	|x| = \langle x, x \rangle^{1/2}.
\end{equation}
The Cauchy--Schwarz inequality states that
\begin{equation}
	|\langle x, y \rangle| \le |x| \, |y|
\end{equation}
for all $x, y \in {\bf R}^n$, and the triangle inequality states
that
\begin{equation}
	|x + y| \le |x| + |y|.
\end{equation}
The standard Euclidean distance on ${\bf R}^n$ is defined by
\begin{equation}
	d(x, y) = |x - y|.
\end{equation}

	Let $C({\bf R}^n)$ denote the vector space of complex-valued
continuous functions on ${\bf R}^n$.  As usual, sums and products of
continuous functions are continuous, and polynomials are continuous.
If $\{f_l\}_{l=1}^\infty$ is a sequence of continuous functions on
${\bf R}^n$ and $f$ is another function on ${\bf R}^n$, then we say
that $\{f_l\}_{l=1}^\infty$ converges to $f$ uniformly on compact sets
if for each compact subset $K$ of ${\bf R}^n$ and each $\epsilon > 0$
there is a positive integer $L$ such that
\begin{equation}
	|f_l(x) - f(x)| < \epsilon
\end{equation}
for all $x \in K$ and $l \ge L$.  It follows from standard results in
analysis that $f$ is also a continuous function in this case.  If two
sequences of continuous functions on ${\bf R}^n$ converge uniformly on
compact subsets of ${\bf R}^n$, then the sums and products of the
functions in the two sequences also converge uniformly on compact
subsets, to the sum and product of the limits of the original
sequences, respectively.

	Let $f$ be a continuous function on ${\bf R}^n$.  The
\emph{support} of $f$ is denoted $\supp f$ and defined to be the
closure of the set of $x \in {\bf R}^n$ such that $f(x) \ne 0$.  This
is equivalent to saying that the complement of the support of $f$ in
${\bf R}^n$ consists of the points $x \in {\bf R}^n$ such that $f$
vanishes on a heighborhood of $x$.  We write $C_{00}({\bf R}^n)$ for
the vector space of continuous functions on ${\bf R}^n$ with compact
support, which is a vector subspace of $C({\bf R}^n)$.  

	A continuous function $f$ on ${\bf R}^n$ is said to vanish at
infinity if for every $\epsilon > 0$ there is a compact subset $K$ of
${\bf R}^n$ such that
\begin{equation}
	|f(x)| < \epsilon
\end{equation}
for all $x \in {\bf R}^n \backslash K$.  The vector space of
continuous functions on ${\bf R}^n$ which vanish at infinity is
denoted $C_0({\bf R}^n)$.  A continuous function $f$ on ${\bf R}^n$
is said to be bounded if there is a nonnegative real number $M$
such that
\begin{equation}
	|f(x)| \le M
\end{equation}
for all $x \in M$.  The vector space of bounded continuous functions
on ${\bf R}^n$ is denoted $C_b({\bf R}^n)$.  Observe that the product
of a continuous function on ${\bf R}^n$ which vanishes at infinity and
a bounded continuous function on ${\bf R}^n$ also vanishes at
infinity.

	Let us say that a continuous function $f$ on ${\bf R}^n$ has
at most polynomial growth if there is a nonnegative real number $C$
and a nonnegative integer $k$ such that
\begin{equation}
	|f(x)| \le C \, (1 + |x|^k)
\end{equation}
for all $x \in {\bf R}^n$.  The vector space of continuous functions
with at most polynomial growth is denoted $C_p({\bf R}^n)$.
Thus we have the inclusions
\begin{equation}
	C_{00}({\bf R}^n) \subseteq C_0({\bf R}^n)
		\subseteq C_b({\bf R}^n) \subseteq C_p({\bf R}^n)
			\subseteq C({\bf R}^n).
\end{equation}

	For $f \in C_b({\bf R}^n)$ we can define the supremum norm
\begin{equation}
	\|f\| = \sup \{|f(x)| : x \in {\bf R}^n\},
\end{equation}
and this satisfies
\begin{equation}
	\|f_1 + f_2\| \le \|f_1\| + \|f_2\|,
		\quad \|f_1 \, f_2\| \le \|f_1\| \, \|f_2\|
\end{equation}
for all $f_1, f_2 \in C_b({\bf R}^n)$.  The supremum metric on
$C_b({\bf R}^n)$ is defined by
\begin{equation}
	\sigma (f_1, f_2) = \|f_1 - f_2\|.
\end{equation}
Convergence of a sequence of functions in $C_b({\bf R}^n)$ with
respect to the supremum metric is equivalent to classical uniform
convergence on all of ${\bf R}^n$.  It is well known that $C_b({\bf
R}^n)$ is complete as a metric space when equipped with the supremum
norm, which is to say that every Cauchy sequence in $C_b({\bf R}^n)$
with respect to the supremum norm converges in $C_b({\bf R}^n)$.

	It is also well known that if a sequence of continuous
functions on ${\bf R}^n$ which vanish at infinity converges uniformly
on ${\bf R}^n$, then the limit also vanishes at infinity.  In other
words, $C_0({\bf R}^n)$ is a closed linear subspace of $C_b({\bf
R}^n)$ with respect to the supremum metric.  Moreover, $C_{00}({\bf
R}^n)$ is dense in $C_0({\bf R}^n)$ with respect to the supremum
metric, which is equivalent to saying that if $f$ is a continuous
function on ${\bf R}^n$ which vanishes at infinity, then there is a
sequence of continuous functions on ${\bf R}^n$ with compact support
which converge to $f$ uniformly.  These approximations of $f$ can be
obtained by multiplying $f$ by continuous functions which are equal to
$1$ on large compact sets, have compact support, and have supremum
norm equal to $1$.  

	Let us say that a sequence of functions $\{f_l\}_{l=1}^\infty$
in $C_b({\bf R}^n)$ converges to another function $f$ in $C_b({\bf
R}^n)$ in the restricted sense if the supremum norms of the $f_l$'s
are all bounded, and if the $f_l$'s converge to $f$ uniformly on
compact subsets of ${\bf R}^n$.  Of course this holds if the $f_l$'s
converge to $f$ uniformly on all of ${\bf R}^n$.  Restricted
convergence is more general, however, and indeed if $f$ is any bounded
continuous function on ${\bf R}^n$, then there is a sequence of
continuous functions on ${\bf R}^n$ with compact support which
converges to $f$ in the restricted sense.  Similarly, if $f$ is any
continuous function on ${\bf R}^n$, then there is a sequence of
continuous functions on ${\bf R}^n$ with compact supports which
converges to $f$ uniformly on compact subsets.  These approximations
of $f$ can be obtained by multiplying $f$ by the same kind of
functions with compact support as in the previous paragraph.

	There is a natural kind of restricted convergence in $C_p({\bf
R}^n)$, which is to say that a sequence $\{f_l\}_{l=1}^\infty$ of
functions in $C_p({\bf R}^n)$ converges to a function $f$ in this
sense if there exist $C, k \ge 0$ such that
\begin{equation}
	|f_l(x)| \le C \, (1 + |x|^k)
\end{equation}
for all $x \in {\bf R}^n$ and $l \ge 1$, and if $\{f_l\}_{l=1}^\infty$
converges to $f$ uniformly on compact sets.  This is equivalent to
saying that there is a $k \ge 0$ such that
\begin{equation}
	\frac{f_l - f}{1 + |x|^{k+1}}
\end{equation}
converges uniformly to $0$ on ${\bf R}^n$ as $l \to \infty$.  Of
course, if a sequence of functions $\{f_l\}_{l=1}^\infty$ in $C_p({\bf
R}^n)$ converges to a function $f$ on ${\bf R}^n$ in this restricted
sense, then $f$ is also in $C_p({\bf R}^n)$.  Every function in
$C_p({\bf R}^n)$ is the limit of a sequence of compactly supported
continuous functions in this restricted sense.

	Let us mention also a notion of restricted convergence in
$C_{00}({\bf R}^n)$, which is that a sequence $\{f_l\}_{l=1}^\infty$
of functions in $C_{00}({\bf R}^n)$ converges to a function $f$ on
${\bf R}^n$ in this sense if there is a single compact subset $K$ of
${\bf R}^n$ such that the support of $f_l$ is contained in $K$ for all
$l$, and if $\{f_l\}_{l=1}^\infty$ converges to $f$ uniformly.
In this case the support of $f$ is also contained in $K$,
so that $f \in C_{00}({\bf R}^n)$.  Under these conditions
uniform convergence is equivalent to uniform convergence on $K$,
with $f$ set to $0$ on ${\bf R}^n \backslash K$.

	Let us write $\mathcal{M}({\bf R}^n)$ for the vector space of
complex-linear functionals on $C_{00}({\bf R}^n)$ which are bounded on
compact subsets of ${\bf R}^n$.  More precisely, an element $\lambda$
of $\mathcal{M}({\bf R}^n)$ is a linear mapping from $C_{00}({\bf
R}^n)$ into the complex numbers with the property that for each
compact subset $K$ of ${\bf R}^n$, there is a nonnegative real number
$L(K)$ such that
\begin{equation}
	|\lambda(f)| \le L(K) \, \|f\|
\end{equation}
whenever $f \in C_{00}({\bf R}^n)$ has support contained in $K$.
Informally we might refer to the elements of $\mathcal{M}({\bf R}^n)$
as measures on ${\bf R}^n$.  We shall write $\mathcal{M}_r({\bf R}^n)$
for the real measures on ${\bf R}^n$, which is to say the $\lambda \in
\mathcal{M}({\bf R}^n)$ such that $\lambda(f)$ is a real number
whenever $f$ is a real-valued continuous function on ${\bf R}^n$.

	If $\lambda \in \mathcal{M}({\bf R}^n)$, then $\lambda$ is
a continuous linear functional on $C_{00}({\bf R}^n)$ with respect to
the kind of restricted convergence defined before.  Namely,
if $K$ is a compact subset of ${\bf R}^n$, $\{f_l\}_{l=1}^\infty$
is a sequence of continuous functions on ${\bf R}^n$ with supports
contained in $K$, and if $f$ is a continuous function on ${\bf R}^n$
with support contained in $K$ such that $\{f_l\}_{l=1}^\infty$
converges uniformly to $f$, then
\begin{equation}
	\lim_{l \to \infty} \lambda(f_l) = \lambda(f).
\end{equation}
This is easy to derive from the boundedness of $\lambda$ on compact
subsets of ${\bf R}^n$.

	Let $\mathcal{M}({\bf R}^n)$ denote the set of nonnegative
linear functionals on $C_{00}({\bf R}^n)$, which is to say the linear
mappings $\lambda$ from $C_{00}({\bf R}^n)$ into the complex numbers
such that for each nonnegative real-valued function $f$ in
$C_{00}({\bf R}^n)$ we have that $\lambda(f)$ is a real number and
\begin{equation}
	\lambda(f) \ge 0.
\end{equation}
If $\lambda \in \mathcal{M}_+({\bf R}^n)$, then $\lambda(f)$ is a real
number for every real-valued function $f$ in $C_{00}({\bf R}^n)$,
because every such function can be written as $f_1 - f_2$ where $f_1$,
$f_2$ are nonnegative real-valued functions in $C_{00}({\bf R}^n)$.

	Suppose that $\lambda \in \mathcal{M}_+({\bf R}^n)$ and
that $\phi$, $f$ are functions in $C_{00}({\bf R}^n)$ such that
$\phi$ is real-valued and
\begin{equation}
	|f(x)| \le \phi(x)
\end{equation}
for all $x \in {\bf R}^n$.  This is equivalent to saying that
\begin{equation}
	\re (\phi(x) - \theta \, f(x)) \ge 0
\end{equation}
for all complex numbers $\theta$ with $|\theta| = 1$ and all
$x \in {\bf R}^n$.  In this case we get that
\begin{equation}
	\re (\lambda(\phi) - \theta \, \lambda(f)) \ge 0
\end{equation}
for all complex numbers $\theta$ such that $|\theta| = 1$, and
therefore that
\begin{equation}
	|\lambda(f)| \le \lambda(\phi)
\end{equation}
in this case.

	In particular it follows that $\lambda$ is bounded on compact
subsets of ${\bf R}^n$, and hence that $\lambda \in \mathcal{M}_r({\bf
R}^n)$.  Indeed, if $K$ is a compact subset of ${\bf R}^n$ and
$\lambda$ is a nonnegative real-valued function in $C_{00}({\bf R}^n)$
such that $\phi(x) = 1$ for all $x \in K$, then one can check that
\begin{equation}
	|\lambda(f)| \le \lambda(\phi) \, \|f\|
\end{equation}
for all $f \in C_{00}({\bf R}^n)$ with support contained in $K$.
Notice that $\mathcal{M}_r({\bf R}^n)$ is a real vector space, and
that $\mathcal{M}_+({\bf R}^n)$ is a cone in $\mathcal{M}({\bf R}^n)$,
in the sense that if $\lambda_1, \lambda_2 \in \mathcal{M}_+({\bf
R}^n)$ and $a_1$, $a_2$ are nonnegative real numbers, then
$a_1 \, \lambda_1 + a_2 \lambda_2 \in \mathcal{M}_+({\bf R}^n)$.
The elements of $\mathcal{M}({\bf R}^n)$ can be described informally
as nonnegative measures on ${\bf R}^n$.

	Suppose that $U$ is an open subset of ${\bf R}^n$, and that
$\lambda \in \mathcal{M}({\bf R}^n)$.  We say that $\lambda$ vanishes
on $U$ if $\lambda(f) = 0$ whenever $f \in C_{00}({\bf R}^n)$ has
support contained in $U$.  If $U_1$, $U_2$ are open subsets of ${\bf
R}^n$, $\lambda \in \mathcal{M}({\bf R}^n)$, and $\lambda$ vanishes on
$U_1$, $U_2$, then $\lambda$ vanishes on $U_1 \cup U_2$.  For if $f
\in C_{00}({\bf R}^n)$ has support contained in $U_1 \cup U_2$, then
$f$ can be written as $f_1 + f_2$ where $f1, f_2 \in C_{00}({\bf
R}^n)$ and the supports of $f_1$, $f_2$ are contained in $U_1$, $U_2$,
respectively, and therefore
\begin{equation}
	\lambda(f) = \lambda(f_1) + \lambda(f_2) = 0.
\end{equation}

	More generally, if $\{U_\alpha\}_{\alpha \in A}$ is a family
of open subsets of ${\bf R}^n$, $\lambda \in \mathcal{M}({\bf R}^n)$,
and $\lambda$ vanishes on each $U_\alpha$, then $\lambda$ vanishes on
the union of the $U_\alpha$'s.  Indeed, let $f$ be a function in
$C_{00}({\bf R}^n)$ whose support is contained in the union of the
$U_\alpha$'s.  Because the support of $f$ is compact by assumption, it
follows that the support of $f$ is contained in the union of finitely
many $U_\alpha$'s.  The result in the previous paragraph may be
applied to obtain that $\lambda(f) = 0$.

	Let $\lambda \in \mathcal{M}({\bf R}^n)$ be given.  The
support of $\lambda$ is defined to be the set of points $x$ in ${\bf
R}^n$ such that $\lambda$ does not vanish on a neighborhood of $x$.
In other words, the complement of the support of $\lambda$ is the set
of points $y \in {\bf R}^n$ such that $\lambda$ vanishes on a
neighborhood of $y$.  This complementary set is open by construction,
and hence the support of $\lambda$ is automatically a closed subset of
${\bf R}^n$.  Also, $\lambda$ vanishes on the complement of its
support in ${\bf R}^n$, by the result mentioned in the previous
paragraph.

	If $\lambda \in \mathcal{M}({\bf R}^n)$ and $\phi \in C({\bf
R}^n)$, then we can define a new measure $\lambda_\phi \in
\mathcal{M}({\bf R}^n)$ by
\begin{equation}
	\lambda_\phi(f) = \lambda(\phi \, f)
\end{equation}
for $f \in C_{00}({\bf R}^n)$.  The support of $\lambda_\phi$ is
automatically contained in the support of $\phi$ as a continuous
function on ${\bf R}^n$.  A measure $\lambda \in \mathcal{M}({\bf
R}^n)$ has compact support if and only
\begin{equation}
	\lambda = \lambda_\phi
\end{equation}
for some $\phi \in C_{00}({\bf R}^n)$.

	If $\phi \in C_{00}({\bf R}^n)$ and $\lambda \in
\mathcal{M}({\bf R}^n)$, then we extend $\lambda_\phi$ to a linear
functional on all of $C({\bf R}^n)$ in an obvious manner.
Specifically, we put
\begin{equation}
	\lambda_\phi(f) = \lambda(\phi \, f)
\end{equation}
for all $f \in C({\bf R}^n)$, and this makes sense because $\phi \, f
\in C_{00}({\bf R}^n)$ since $\phi \in C_{00}({\bf R}^n)$.  This
extension to $C({\bf R}^n)$ is continuous in the sense that if
$\{f_l\}_{l=1}^\infty$ is a sequence of continuous functions on ${\bf
R}^n$ which converges uniformly to the continuous function $f$ on
compact subsets of ${\bf R}^n$, then
\begin{equation}
	\lim_{l \to \infty} \lambda_\phi(f_l) = \lambda_\phi(f),
\end{equation}
essentially because we only need the uniform convergence of the
$f_l$'s to $f$ on the support of $\phi$.  This continuity property
uniquely characterizes the extension of $\lambda_\phi$ to $C({\bf
R}^n)$, since $C_{00}({\bf R}^n)$ is dense in $C({\bf R}^n)$ with
respect to uniform convergence on compact subsets of ${\bf R}^n$.

	A measure $\lambda \in \mathcal{M}({\bf R}^n)$ is said to
be bounded if there is a nonnegative real number $L$ such that
\begin{equation}
	|\lambda(f)| \le L \, \|f\|
\end{equation}
for all $f \in C_{00}({\bf R}^n)$.  We write $\mathcal{M}_b({\bf
R}^n)$ for the vector subspace of $\mathcal{M}({\bf R}^n)$ of bounded
measures.  For $\lambda \in \mathcal{M}_b({\bf R}^n)$ we define
the norm of $\lambda$ by
\begin{equation}
	\|\lambda\|_* = \sup \{|\lambda(f)| : f \in C_{00}({\bf R}^n),
						\|f\| \le 1\}.
\end{equation}
This is the same as saying that $\|\lambda\|_*$ is the smallest choice
of $L$ for which the previous inequality holds.

	It is easy to see that if $\lambda \in \mathcal{M}_b({\bf
R}^n)$ and $a$ is a complex number, then
\begin{equation}
	\|a \, \lambda\|_* = |a| \, \|\lambda\|_*.
\end{equation}
If $\lambda_1, \lambda_2 \in \mathcal{M}_b({\bf R}^n)$, then
\begin{equation}
	\|\lambda_1 + \lambda_2\|_* \le \|\lambda_1\|_* + \|\lambda_2\|_*.
\end{equation}
If $\lambda \in \mathcal{M}_b({\bf R}^n)$ and $\phi \in C_b({\bf
R}^n)$, then it is easy to check that $\lambda_\phi \in
\mathcal{M}_b({\bf R}^n)$ and that
\begin{equation}
	\|\lambda_\phi\|_* \le \|\phi\| \, \|\lambda\|_*.
\end{equation}
This inequality is dual to the one that says that $\|\phi \, f\| \le
\|\phi\| \, \|f\|$ when $\phi, f \in C_b({\bf R}^n)$.

	Now suppose that $\lambda \in \mathcal{M}_b({\bf R}^n)$ and
that $\phi_1, \ldots, \phi_l \in C_b({\bf R}^n)$, and let us show that
\begin{equation}
	\sum_{j=1}^l \|\lambda_{\phi_j}\|_*
		\le \biggl\|\sum_{j=1}^l |\phi_j| \biggr\|
				\, \|\lambda\|_*.
\end{equation}
This is equivalent to saying that if $f_1, \ldots, f_l \in C_{00}({\bf
R}^n)$ and $\|f_j\| \le 1$ for $j = 1, \ldots, l$, then
\begin{equation}
	\sum_{j=1}^l |\lambda (\phi_j \, f_j)|
		\le \biggl\|\sum_{j=1}^l |\phi_j| \biggr\| \, \|\lambda\|_*.
\end{equation}
This is equivalent in turn to saying that for such $f_1, \ldots, f_l$
and for complex numbers $\theta_1, \ldots, \theta_l$ with $|\theta_j|
= 1$ for $j = 1, \ldots, l$ we have that
\begin{equation}
	\biggl|\sum_{j=1}^l \theta_j \lambda (\phi_j \, f_j) \biggr|
		\le \biggl\|\sum_{j=1}^l |\phi_j| \biggr\| \, \|\lambda\|_*.
\end{equation}
Since $\lambda$ is linear, this reduces to
\begin{equation}
  \biggl|\lambda \biggl(\sum_{j=1}^n \theta_j \, \phi_j \, f_j \biggr)\biggr|
	\le \biggl\|\sum_{j=1}^l |\phi_j| \biggr\| \, \|\lambda\|_*,
\end{equation}
which holds because
\begin{equation}
	\biggl\|\sum_{j=1}^l \theta_j \, \phi_j \, f_j \biggr\|
		\le \biggl\|\sum_{j=1}^l |\phi_j| \biggr\|
\end{equation}
under our conditions on the $f_j$'s and the $\theta_j$'s.

	Let $\lambda$ be a bounded measure on ${\bf R}^n$, and
let $\epsilon > 0$ be given.  Let $f$ be a continuous function
on ${\bf R}^n$ with compact support such that $\|f\| \le 1$ and
\begin{equation}
	|\lambda(f)| \ge \|\lambda\|_* - \epsilon.
\end{equation}
Let $\phi$ be a real-valued continuous function on ${\bf R}^n$ with
compact support such that $0 \le \phi(x) \le 1$ for all $x \in {\bf
R}^n$ and $\phi(x) \, f(x) = f(x)$ for all $x \in {\bf R}^n$.
Thus
\begin{equation}
	\lambda_\phi(f) = \lambda(f)
\end{equation}
and therefore
\begin{equation}
	\|\lambda_\phi\|_* \ge |\lambda_\phi(f)| \ge \|\lambda\|_* - \epsilon.
\end{equation}
We also have that $\lambda_\phi + \lambda_{1 - \phi} = \lambda$ and
\begin{equation}
	\|\lambda\|_* = \|\lambda_\phi\|_* + \|\lambda_{1 - \phi}\|_*
\end{equation}
by the result of the preceding paragraph, and therefore
\begin{equation}
	\|\lambda - \lambda_\phi\|_* = \|\lambda_{1 - \phi}\|_* \le \epsilon.
\end{equation}
In other words, bounded measures on ${\bf R}^n$ can be approximated
by measures with compact support in the dual norm on measures.

	Suppose that $\lambda$ is a bounded measure on ${\bf R}^n$,
and that $\{\lambda_j\}_{j=1}^\infty$ is a sequence of measures on
${\bf R}^n$ with compact support which converge to $\lambda$ in the
dual norm, which is to say that 
\begin{equation}
	\lim_{j \to \infty} \|\lambda_j - \lambda\|_* = 0.
\end{equation}
Let $f$ be a bounded continuous function on ${\bf R}^n$.  Thus
$\lambda_j(f)$ is defined for all $j$, since the $\lambda_j$'s have
compact support.  Moreover,
\begin{equation}
	|\lambda_j(f) - \lambda_l(f)| 
		\le \|\lambda_j - \lambda_l\|_* \, \|f\|
\end{equation}
for all positive integers $j$, $l$, and it follows that
$\{\lambda_j(f)\}_{j=1}^\infty$ is a Cauchy sequence of complex
numbers.  Every Cauchy sequence of complex numbers converges,
and so we may define $\lambda(f)$ for $f \in C_b({\bf R}^n)$ by
\begin{equation}
	\lambda(f) = \lim_{j \to \infty} \lambda_j(f).
\end{equation}
Of course this agrees with the initial definition of $\lambda(f)$
when $f \in C_{00}({\bf R}^n)$.

	If we think of $C_0({\bf R}^n)$ as a metric space equipped
with the supremum metric, then a bounded measure on ${\bf R}^n$ is a
linear functional on the dense subspace $C_{00}({\bf R}^n)$ of
$C_0({\bf R}^n)$ which is continuous with respect to the supremum
metric, and in fact uniformly continuous, because of linearity.
By standard results about metric spaces there is a unique continuous
extension of the linear functional to $C_0({\bf R}^n)$, and the
extension is linear and bounded with the same norm as the original.
The extension described in the previous paragraph agrees with this
one on $C_0({\bf R}^n)$, and goes further.

	More precisely, the extension of a bounded measure $\lambda$
on ${\bf R}^n$ to $C_b({\bf R}^n)$ described above enjoys a stronger
continuity property, which is continuity with respect to restricted
convergence in $C_b({\bf R}^n)$.  This means that if
$\{f_m\}_{m=1}^\infty$ is a sequence of bounded continuous functions
on ${\bf R}^n$ which are uniformly bounded on ${\bf R}^n$ and which
converge uniformly on compact subsets of ${\bf R}^n$ to a function $f
\in C_b({\bf R}^n)$, then the extension of $\lambda$ satisfies
\begin{equation}
	\lim_{m \to \infty} \lambda (f_m) = \lambda (f).
\end{equation}
One can check this using the approximation of $\lambda$ in the
dual norm by measures on ${\bf R}^n$ with compact support.  Of course
for a measure with compact support one has continuity with respect to
uniform convergence on compact subsets of ${\bf R}^n$.

	Because $C_{00}({\bf R}^n)$ is dense in $C_b({\bf R}^n)$ with
respect to this kind of restricted convergence, we have that our
extension of $\lambda$ to $C_b({\bf R}^n)$ is uniquely determined by
the initial definition of $\lambda$ on $C_{00}({\bf R}^n)$ and this
continuity with respect to restricted convergence in $C_b({\bf R}^n)$.
In particular the extension of $\lambda$ to $C_b({\bf R}^n)$ does not
depend on the choice of sequence of measures on ${\bf R}^n$ with
compact support which approximate $\lambda$, although this is easy to
check directly too.  Of course this extension of $\lambda$ to
$C_b({\bf R}^n)$ is linear and satisfies
\begin{equation}
	|\lambda(f)| \le \|\lambda\|_* \, \|f\|
\end{equation}
for all $f \in C_b({\bf R}^n)$.

\section{Fourier transforms}
\label{Fourier transforms}
\setcounter{equation}{0}

	Let $\lambda$ be a bounded measure on ${\bf R}^n$.  For
each $\xi \in {\bf R}^n$, let $e_\xi$ denote the bounded continuous
function on ${\bf R}^n$ given by
\begin{equation}
	e_\xi(x) = \exp (- 2 \pi i \xi \cdot x),
\end{equation}
with
\begin{equation}
	\xi \cdot x = \sum_{j=1}^n \xi_j \, x_j,
\end{equation}
$x = (x_1, \ldots, x_n)$, $\xi = (\xi_1, \ldots, \xi_n)$.  As in the
previous section, $\lambda(f)$ is defined when $f$ is a bounded
continuous function on ${\bf R}^n$, and we define the Fourier
transform of $\lambda$ by
\begin{equation}
	\widehat{\lambda}(\xi) = \lambda(e_\xi).
\end{equation}

	If $\xi \in {\bf C}^n$ we can define $e_\xi(x)$ as a
continuous function on ${\bf R}^n$ in exactly the same manner as
above.  If $\lambda$ is a measure on ${\bf R}^n$ with compact support,
then $\lambda(h)$ is defined for all continuous functions on ${\bf
R}^n$.  In particular, $\widehat{\lambda}(\xi) = \lambda(e_\xi)$ is
then defined for all $\xi \in {\bf C}^n$, and in fact it is a complex
analytic function of $\xi$, which is to say that
$\widehat{\lambda}(\xi)$ can be expressed as a power series in the
$\xi_j$'s.  In general for bounded measures $\lambda$ it may be
possible to define the Fourier transform of $\lambda$ for some or all
$\xi \in {\bf C}^n$, depending on the behavior of $\lambda$.

	If $\lambda$ is a bounded measure on ${\bf R}^n$, then we have
that
\begin{equation}
	|\widehat{\lambda}(\xi)| \le \|\lambda\|_*
		\quad\hbox{for all } \xi \in {\bf R}^n.
\end{equation}
Moreover, $\widehat{\lambda}(\xi)$ is a uniformly continuous function
on ${\bf R}^n$.  Indeed, if $\lambda$ has compact support in ${\bf
R}^n$, then its Fourier transform is quite a bit more regular than
that.  In general, a bounded measure can be approximated by measures
with compact support in the dual norm, and this means that the Fourier
transform of a bounded measure can be approximated in the supremum
norm by the Fourier transforms of measures with compact support, so
that the uniform continuity of the Fourier transform of a bounded
measure follows from the regularity of the Fourier transforms of
measures with compact support.

	For each continuous function $f$ on ${\bf R}^n$ and
each $v \in {\bf R}^n$ define $\tau_v(f)$ to be the continuous
function on ${\bf R}^n$ given by
\begin{equation}
	\tau_v(f)(y) = f(v - y).
\end{equation}
If $\lambda$ is a measure on ${\bf R}^n$ and $f$ is a continuous
function on ${\bf R}^n$, we would like to define the convolution
of $\lambda$ and $f$ to be the function on ${\bf R}^n$ expressed
by the formula
\begin{equation}
	(\lambda * f)(v) = \lambda(\tau_v(f)).
\end{equation}
This makes sense if $f$ is a continuous function on ${\bf R}^n$ with
compact support and $\lambda$ is any measure on ${\bf R}^n$, in which
case $\lambda * f$ is a continuous function on ${\bf R}^n$.

	Alternatively, if $\lambda$ is a measure on ${\bf R}^n$ with
compact support, and if $f$ is any continuous function on ${\bf R}^n$,
then $(\lambda * f)(v)$ is again defined.  One can check that $\lambda
* f$ is also continuous in this event, and indeed for a compact set of
$v$'s $(\lambda * f)(v)$ will involve the values of $f$ only on a
compact subset of ${\bf R}^n$.  If both $\lambda$ and $f$ have compact
support, then $\lambda * f$ has compact support as well.

	Now suppose that $\lambda$ is a bounded measure on ${\bf R}^n$,
and that $f$ is a bounded continuous function on ${\bf R}^n$.
In this case $\lambda * f$ can be defined again, and one can check
that $\lambda * f$ is continuous using the fact that $\lambda$
is continuous on $C_b({\bf R}^n)$ with respect to restricted convergence
of sequences of bounded continuous functions.  Of course
\begin{equation}
	\|\lambda * f \| \le \|\lambda\|_* \, \|f\|.
\end{equation}
One can also look at the continuity of $\lambda * f$ in terms
of approximations of $\lambda$ in the dual norm by measures with
compact support, which lead to approximations of $\lambda * f$
by convolutions of $f$ with measures with compact support in the
supremum metric on $C_b({\bf R}^n)$.

	If $\lambda$ is a bounded measure on ${\bf R}^n$ and
$\xi \in {\bf R}^n$, then we have that
\begin{equation}
	\lambda * e_\xi = \widehat{\lambda}(-\xi) \, e_\xi.
\end{equation}
Thus the exponential function $e_xi$ is an eigenfunction for the
linear operator
\begin{equation}
	T_\lambda(f) = \lambda * f
\end{equation}
on $C_b({\bf R}^n)$, and the corresponding eigenvalue is given by the
Fourier transform of $\lambda$ at $-\xi$.  If $\lambda$ has compact
support, then we can think of $T_\lambda$ as a linear mapping from
continuous functions on ${\bf R}^n$ to themselves, and $e_\xi$
is an eigenfunction of $T_\lambda$ for all $\xi \in {\bf C}^n$,
with eigenvalue equal to the Fourier transform of $\lambda$
at $-\xi$.

\end{document}